\documentclass[english]{article}
\usepackage[T1]{fontenc}
\usepackage[latin1]{inputenc}
\usepackage{slashed}
\usepackage{amssymb}

\newcommand{\Dh}{Dirac-harmonic map}
\usepackage{babel}
\makeatother \evensidemargin46mm \oddsidemargin10mm
\title{Some explicit constructions of Dirac-harmonic maps
} \footnotesize{\author{J\"urgen  Jost, Xiaohuan
Mo\thanks{Supported by the National Natural Science Foundation of
China 10771004}\, and Miaomiao Zhu\thanks{Supported by IMPRS
``Mathematics in
the Sciences'' and the Klaus Tschira Foundation}\\
Max Planck Institute for Mathematics in the Sciences,\\
Inselstr. 22, D-04103 Leipzig, Germany\\
E-mail:  jjost@mis.mpg.de, miaomiao.zhu@mis.mpg.de\\
Key Laboratory of Pure and Applied Mathematics\\
School of Mathematical Sciences\\ Peking  University,
Beijing 100871, China\\
E-mail:  moxh@pku.edu.cn}}

\begin{document}
\maketitle
\begin{center}{\bf Abstract} \end{center}
\begin{center}
\begin{minipage}{125mm}
We  construct explicit examples of Dirac-harmonic maps
$(\phi,\,\psi)$ between Riemannian manifolds $(M,g)$ and $(N,g')$
which are non-trivial in the sense that $\phi$ is not harmonic.
When $\dim M=2$, we also produce examples where $\phi$ is
harmonic, but not conformal, and $\psi$ is non-trivial.

{\bf Key words and phrases:} Dirac-harmonic map, twistor spinor, totally umbilical

{\bf 1991 Mathematics Subject Classification:} 58E20.

\end{minipage}
\end{center}

\section{Introduction}

A Dirac-harmonic map is a pair that  couples a map between
Riemannian manifolds with a nonlinear spinor field along that map
[10]. \Dh s   arise from the supersymmetric nonlinear sigma model of
quantum field theory [12]. They are a generalization and combination
of harmonic maps and harmonic spinors while preserving
the essential properties of the former.\\
Both harmonic maps and harmonic spinors have been extensively studied. See, for instance [13,\,15]. In particular, many non-trivial examples of harmonic maps and harmonic spinors are known
[2,\,3,\,4,\,13]. A harmonic map and a vanishing spinor, or conversely
a constant map and a harmonic spinor constitute an example of a \Dh
. A natural question then is whether there exist other examples that
couple a map and a spinor in a non-trivial manner.
The purpose of this paper therefore is to manufacture non-trivial examples of Dirac-harmonic maps
between Riemannian manifolds. For hypersurfaces in a  Riemannian manifold of constant
sectional curvature, we prove the following:

\

{\bf Theorem 1}\,\,{\em Let $M$ be an $n$-dimensional manifold
which is immersed in an $(n+1)$-dimensional Riemannian manifold
$N(c)$ of constant sectional curvature $c$. Assume $\Phi$ is a
harmonic spinor on $M$, and $\Psi\in\Gamma(\Sigma M)$ satisfies
$$
-2c Re\langle\Phi,\,\Psi\rangle\nu=H.
\eqno(1)
$$
where $H$ is the mean curvature vector field of $\phi$, $\nu$ is the unit normal field of $\phi$ and $\Sigma M$ is the spinor bundle of $M$.
We define a spinor field $\psi$ along the immersion $\phi$ by
$$
\psi=\Sigma_{\alpha} \epsilon_{\alpha}\cdot\Psi\otimes
\phi_*(\epsilon_{\alpha})+\Phi\otimes\nu
$$
where $\epsilon_{\alpha}$ is a local
orthonormal basis of $M$.

(i) If $n=2$,  $\phi$ is minimal and $\Psi$ satisfies
$$
\epsilon_{1}\cdot\nabla_{\epsilon_{1}}\Psi-\epsilon_{2}\cdot\nabla_{\epsilon_{2}}\Psi=\lambda_1\Phi
\eqno(2)
$$
where $\lambda_{1}$ is the principal curvature in the principal direction $\epsilon_{1}$, then $(\phi,\,\psi)$ is Dirac-harmonic.

(ii) If $n\geq 3$,  $\phi$ is totally umbilical and $\Psi$ is a twistor spinor  satisfying
$$
\slashed{\partial}\Psi=-\frac{n\langle H,\,\nu\rangle}{n-2}\Phi
\eqno(3)
$$
then $(\phi,\,\psi)$ is Dirac-harmonic.

}

\

Using Theorem 1, we can construct many Dirac-harmonic maps
$(\phi,\,\psi)$ from $\mathbb{R}^n$ into $\mathbb{H}^{n+1}(-1)$
where $n\geq 3$ and $\phi:\mathbb{R}^n\to\mathbb{H}^{n+1}(-1)$ is
not harmonic (see Section 5, Example 3).

Finding explicit non-trivial explicit solutions of (2) and (3) turns out to be
difficult. However in some     special cases, we are able to get
the non-trivial solutions, as in Example 3.

Let us take a look at the following special case of (i) of Theorem 1: when $\Phi=0$, then
$(\phi,\,\psi)$ is Dirac-harmonic if $\phi:M\to N(c)$ is minimal and $\Psi$
is a twistor spinor. In fact, in the general case, we have the following:

\

{\bf Theorem 2}\,\,{\em  Let $M$ be a Riemann
surface and $N$ a Riemannian manifold. Assume $\phi:M\to N$
is a harmonic map and $\Psi\in
\Gamma(\Sigma M)$ is a twistor spinor. We define a spinor field $\psi_{\phi,\,\Psi}$ along map $\phi$ by
$$
\psi_{\phi,\,\Psi}:=\Sigma_{\alpha} \epsilon_{\alpha}\cdot\Psi\otimes
\phi_*(\epsilon_{\alpha})\eqno(4)
$$
where $\epsilon_{\alpha}$ ($\alpha=1,\,2$) is a local
orthonormal basis of  $\,\,M$. Then $(\phi,\,\psi_{\phi,\,\Psi})$ is
a Dirac-harmonic map.}

\

By using Theorem 2, we can manufacture Dirac-harmonic maps
$(\phi,\,\psi_{\phi,\,\Psi})$ from a surface for a (not necessarily
conformal) map $\phi$ (see Section 4). Theorem 2 generalizes the
result of [10] that was derived for the special case when both
source and target manifolds are two-dimensional spheres.

Finally, by investigating spinor fields along a hypersurface with
two constant principal curvatures in a Riemannian manifold of
constant curvature, we get Dirac-harmonic maps $(\phi,\,\psi)$
from surfaces for which $\phi$ is not harmonic (see Section 6 and
Section 7).

Let us describe our construction. Let $M:=S^1(r)\times
H^{1}(\sqrt{R^2+r^2})$ be a hyperbolic  surface of revolution (see
Section 6 for definitions). Let $a$ and $b$ be arbitrary complex
constants and $m$ be an arbitrary non-negative integer. For each
$k\in \{0,\,\pm 1\,\cdots,\pm m\}$, let $c_k$ and $d_k$ be complex
constants satisfying
$$
Re(a\bar{d_0}+\bar{b}{c_0})=\frac{\sqrt{R^2+r^2}(R^2+2r^2)}{2rR}
\eqno(5)
$$
and
$$
a\bar{d_k}+\bar{b}c_{-k}=0. \eqno(6)
$$
We obtain the following result (see Section 7):

\

{\bf Theorem 3}\,\,{\em Let $\phi:M\hookrightarrow H^{3}(R)$ be an
isometric immersion from $M$ into a hyperbolic space and
$\psi\in\Gamma(\Sigma M\otimes \phi^{-1}TH^{3}(R))$ defined by
$$
\psi= \epsilon_{1}\cdot\Psi\otimes
\phi_*(\epsilon_{1})-\frac{r^2}{R^2+r^2}\epsilon_{2}\cdot\Psi\otimes
\phi_*(\epsilon_{2})+\chi\otimes\nu
$$
where
$$
\nu(\theta,\,t)=-\left(\frac{\sqrt{R^2+r^2}}R\cos\frac{\theta}r,\,\frac{\sqrt{R^2+r^2}}R\sin\frac{\theta}r,\,\frac{r}R\sinh\frac
t{\sqrt{R^2+r^2}},\,\frac{r}R\cosh\frac t{\sqrt{R^2+r^2}}\right)
$$ is a unit normal vector of $M$, and
$\chi
=\left(\begin{array}{l}  a\\
b
\end{array}
\right)$
$$
\Psi(\theta,\,t)=i\frac{\sqrt{R^2+r^2}}{rR} t\left(\begin{array}{l}  b\\
a
\end{array}
\right) +\sum_{k=-m}^m e^{i\frac kr\theta}
\left(\begin{array}{l}  d_ke^{-\frac krt}\\
c_ke^{\frac krt}
\end{array}
\right)
$$
are the spinors on $M$ with respect to the "untwisted" spinor
bundle on $M$ satisfying (5) and (6).
$\{\epsilon_1,\,\epsilon_2\}$ is a local orthonormal basis of $M$
such that $$
\epsilon_1(\theta,\,t)=\left(-\sin\frac{\theta}r,\,\cos\frac{\theta}r,\,0,\,0\right)
$$ is a principal curvature $\frac{\sqrt{R^2+r^2}}{rR}$ direction and $$
\epsilon_2(\theta,\,t)=\left(0,\,0,\,\cosh\frac
t{\sqrt{R^2+r^2}},\,\sinh\frac t{\sqrt{R^2+r^2}}\right)
$$ is  a principal curvature $\frac{r}{R\sqrt{R^2+r^2}}$ direction. Then $(\phi,\,\psi)$ is a Dirac-harmonic
map from $M$ into $H^{3}(R)$ for which $\phi$ is not harmonic. }

\

The proofs of our results are essentially of an algebraic nature.
They carefully match the algebraic structure of the curvature term
in the Dirac-harmonic map equation, as displayed in the next
section, with the special properties of twistor spinors or those of
particular submanifolds defined in terms of ambient curvature
properties in spaces of constant curvature.

\section{Dirac-harmonic maps}

Let $(N,\,h)$ be a Riemannian manifold of dimension $n'$, $(M,\,g)$ be
an $n$-dimensional Riemannian manifold with fixed spin structure,
$\Sigma M$ its spinor bundle, on which we have a Hermitian metric
$\langle\cdot,\,\cdot\rangle$  induced by the Riemannian metric
$g(\cdot,\,\cdot)$ of $M$. Let $\phi$ be a smooth map from $(M,\,g)$ to $(N,\,h)$ and $\phi^{-1}TN$
the pull-back bundle of $TN$ by $\phi$. On the twisted bundle
$\Sigma M\otimes\phi^{-1}TN$ there is a metric (still denoted by $\langle\cdot,\,\cdot\rangle$)
induced from the metrics on $\Sigma M$ and $\phi^{-1}TN$. There is also a natural connection
$\tilde{\nabla}$ on $\Sigma M\otimes\phi^{-1}TN$ induced from those on
$\Sigma M$ and $\phi^{-1}TN$ (which in turn come from the Levi-Civita
connections of $(M,\,g)$ and $(N,\,h)$, resp.).

For $X\in \Gamma(TM)$, $\xi\in \Gamma(\Sigma M)$, denote by $X\cdot\xi$ their Clifford product, which satisfies
the  skew-symmetry relation
$$
\langle X\cdot\xi,\,\eta\rangle=-\langle\xi,\,X\cdot\eta\rangle
\eqno (7)
$$
as well as the Clifford relations
$$
X\cdot Y\cdot\psi+Y\cdot X\cdot\psi=-2g(X,\,Y)\psi
$$
for $X,\,Y\in \Gamma(TM)$, $\xi,\,\eta\in \Gamma(\Sigma M)$.

Let $\psi$ be a section of the bundle $\Sigma M\otimes\phi^{-1}TN$.
The {\em Dirac operator along the map $\phi$} is defined as
$$
\slashed{D}\psi:=\epsilon_{\alpha}\cdot\tilde{\nabla}_{\epsilon_{\alpha}}\psi
$$
where $\epsilon_{\alpha}$ is a
local orthonormal basis of $\,M$.
For more details about the spin bundle and Dirac operator, we refer to
[17, 21].

Set
$$
\chi:=\{(\phi,\,\psi)\,|\,\,\phi\in
C^{\infty}(M,\,N)\,\,\mbox{and}\,\,\psi\in C^{\infty}(\Sigma
M\otimes\phi^{-1}TN)\}.
$$
On $\chi$, we consider the following functional
$$
L(\phi,\,\psi):=\frac 12\int_M\left[|d\phi|^2+\langle\psi,\,\slashed{D}\psi\rangle\right]{}^*1_M.
$$
This functional couples the two fields $\phi$ and $\psi$ because the
operator $\slashed{D}$ depends on the map $\phi$.
The Euler-Lagrange equations of $L(\phi,\,\psi)$ then also couple the
two fields; they are:
$$
\tau(\phi)=\mathcal{R}(\phi,\,\psi)\eqno(8)
$$
and
$$
\slashed{D}\psi=0\eqno(9)
$$
where $\tau(\phi):={\rm trace} \nabla d\phi$ is the tension field
of the map $\phi$ and $\mathcal{R}(\phi,\,\psi)$ is defined by
$$
\mathcal{R}(\phi,\,\psi)=\frac
12R^i{}_{jkl}\langle\psi^k,\,\nabla\phi^j\cdot\psi^l\rangle\frac{\partial}{\partial
y^i},
$$
where
$$
\psi=\psi^i\otimes\frac{\partial}{\partial
y^i},
$$
$$
(d\phi)^{\sharp}=\nabla \phi^i\otimes\frac{\partial}{\partial
y^i},
$$
$$
R^{\phi^{-1}TN}\left(\frac{\partial}{\partial
y^k},\,\frac{\partial}{\partial
y^l}\right)\frac{\partial}{\partial
y^j}=R^i{}_{jkl}\frac{\partial}{\partial
y^i}
$$
where $^{\sharp}:T^*M\otimes\phi^{-1}TN\to TM\otimes\phi^{-1}TN$ is
the standard (``musical'') isomorphism obtained from the Riemannian
metric $g$.

Solutions $(\phi,\,\psi)$ to (8) and (9) are called {\em
Dirac-harmonic maps} from $M$ into $N$ [9].\\

We now start with some differential geometric identities:
Let $\epsilon_{\alpha}$ be a local
orthonormal basis of  $M$.
By using the Clifford relations we have
$$
\epsilon_{\alpha}\cdot
\epsilon_{\beta}\cdot{\psi}=(-1)^{\delta_{\alpha\beta}+1}\epsilon_{\beta}\cdot
\epsilon_{\alpha}\cdot{\psi} = \left\{
\begin{array}{ccl}
-{\psi}, &&\alpha=\beta\\
-\epsilon_{\beta}\cdot \epsilon_{\alpha}\cdot{\psi},&& \alpha\neq\beta
\end{array}
\right.\eqno(10)
$$
for ${\psi}\in \Gamma(\Sigma M)$.

\

{\bf Lemma 2.1}\,{\em $\mathcal{R}(\phi,\,\psi)\in
  \Gamma(\phi^{-1}TN)$}; in particular, it is real.

\

{\bf Proof}\,\, For any (not necessarily orthonormal) frame
$\{\epsilon_i\}$ on $\phi^{-1}TN$, we put
$$
\psi={\psi}^a\otimes\epsilon_a, \eqno(11)
$$
$$
(d\phi)^{\sharp}=\nabla {\phi}^a\otimes\epsilon_a, \eqno(12)
$$
$$
R^{\phi^{-1}TN}(\epsilon_a,\,\epsilon_b)\epsilon_c={R}^d{}_{abc}\epsilon_d
$$
where $^{\sharp}:T^*M\otimes\phi^{-1}TN\to TM\otimes\phi^{-1}TN$
is the musical isomorphism as before. Take
$$
\epsilon_a=u^i_a\frac{\partial}{\partial y^i},
$$
then
$$
\psi^i=u_a^i{\psi}^a,\qquad
\nabla\phi^i=u^i_a\nabla{\phi^a},\qquad u^j_a u^k_b a^l_c
R^i{}_{jkl}={R}^d{}_{abc}u^i_d.
$$
A simple calculation gives following
$$
R^i{}_{jkl}\langle\psi^k,\,\nabla\phi^j\cdot\psi^l\rangle\frac{\partial}{\partial
y^i}={R}^a{}_{bcd}\left(\phi(x)\right)\langle{\psi}^c,\,\nabla{\phi^b}\cdot{{\psi}}^d\rangle\epsilon_a\left(\phi(x)\right).
\eqno(13)
$$
It follows that the definition of $\mathcal{R}(\phi,\,\psi)$ is
independent of the choice of frame. Moreover, from the skew-symmetry
of $R^i{}_{jkl}$ with respect to the induces $k$ and $l$, we have
$$
\begin{array}{ccl}
\overline{\frac
12R^i{}_{jkl}\langle\psi^k,\,\nabla\phi^j\cdot\psi^l\rangle} &=&
\frac 12R^i{}_{jkl}\langle\nabla\phi^j\cdot\psi^l,\,\psi^k\rangle\\
&=& \frac
12R^i{}_{jlk}\langle\nabla\phi^j\cdot\psi^k,\,\psi^l\rangle\\
&=& -\frac
12R^i{}_{jkl}\langle\nabla\phi^j\cdot\psi^k,\,\psi^l\rangle=\frac
12R^i{}_{jkl}\langle\psi^k,\,\nabla\phi^j\cdot\psi^l\rangle.
\end{array}
$$
It follows that $\mathcal{R}(\phi,\,\psi)$ is well-defined vector
field on $\phi^{-1}TN$, i.e., $\mathcal{R}(\phi,\,\psi) \in \Gamma (
\phi^{-1}TN)$. \hfill $\Box$

\

A spinor (field) $\Psi\in\Gamma(\Sigma M)$ is called a {\em twistor
spinor} if $\Psi$ belongs to the kernel of the twistor operator,
equivalently,
$$
\nabla_{X}\Psi+\frac 1nX\cdot\slashed{\partial}\Psi=0\qquad \forall X\in\Gamma(TM)
$$
where $n$ is the dimension of Riemannian manifold $M$, $\Sigma M$ is
the associated spinor bundle of $M$ and $\slashed{\partial}$ is the
usual Dirac operator (cf. [1,\, 14,\, 20,\, 23]).

In fact the concept of a twistor spinor (in particular, a Killing
spinor) is motivated by theories from physics, like
 general relativity, $11$-dimensional (resp. $10$-dimensional)
supergravity theory, supersymmetry
(see, for example [5,\, 8,\, 11]).

\section{Dirac-harmonic maps from surfaces\,\, I}

In this section, we consider {\em two}-dimensional Riemannian
manifolds $(M,\,g)$. Since a metric on a two-dimensional Riemannian
manifold defines a conformal structure, we then also have the
structure of a Riemann surface. In fact, since the functional $L$
and its critical points, the \Dh s are conformally invariant (see
[10]), in our subsequent considerations, we only need the conformal
structure in place of the full Riemannian metric $g$.

\

{\bf Lemma 3.1}\, {\em Let $\Psi$ be a section of $\Sigma M$. Then
$\langle \epsilon_{\alpha}\cdot\Psi,\, \epsilon_{\beta}\cdot
\epsilon_{\gamma}\cdot\Psi\rangle$ is purely imaginary for any
$\alpha,\,\beta,\,\gamma$ where $\epsilon_{\alpha}$ ($\alpha=1,\,2$) is a
local orthonormal basis of $\,\,M$.}

\

{\bf Proof:}\,\,\, For the Hermitian product
$\langle\cdot,\,\cdot\rangle$ on the spinor bundle $\Sigma M$, we
have
$$
\begin{array}{ccl}
\overline{\langle \epsilon_{\alpha}\cdot\Psi,\, \epsilon_{\beta}\cdot
\epsilon_{\gamma}\cdot\Psi\rangle}
&=&
\langle \epsilon_{\beta}\cdot
\epsilon_{\gamma}\cdot\Psi,\, \epsilon_{\alpha}\cdot\Psi\rangle\\
&=& -\langle \epsilon_{\gamma}\cdot\Psi,\, \epsilon_{\beta}\cdot
\epsilon_{\alpha}\cdot\Psi\rangle\\
&=& -(-1)^{\delta_{\alpha\beta}+1}\langle \epsilon_{\gamma}\cdot\Psi,\,
\epsilon_{\alpha}\cdot \epsilon_{\beta}\cdot\Psi\rangle\\
&=& (-1)^{\delta_{\alpha\beta}+1}\langle \epsilon_{\alpha}\cdot
\epsilon_{\gamma}\cdot\Psi,\, \epsilon_{\beta}\cdot\Psi\rangle\\
&=& (-1)^{\delta_{\alpha\beta}+1}(-1)^{\delta_{\gamma\alpha}+1}
\langle \epsilon_{\gamma}\cdot \epsilon_{\alpha}\cdot\Psi,\,
\epsilon_{\beta}\cdot\Psi\rangle\\
&=& (-1)^{\delta_{\alpha\beta}+\delta_{\gamma\alpha}} \langle
\epsilon_{\gamma}\cdot \epsilon_{\alpha}\cdot\Psi,\, \epsilon_{\beta}\cdot\Psi\rangle\\
&=& -(-1)^{\delta_{\alpha\beta}+\delta_{\gamma\alpha}} \langle
\epsilon_{\alpha}\cdot\Psi,\, \epsilon_{\gamma}\cdot \epsilon_{\beta}\cdot\Psi\rangle\\
&=&
-(-1)^{\delta_{\beta\gamma}+1}(-1)^{\delta_{\alpha\beta}+\delta_{\gamma\alpha}}
\langle \epsilon_{\alpha}\cdot\Psi,\, \epsilon_{\beta}\cdot
\epsilon_{\gamma}\cdot\Psi\rangle\\
&=&
(-1)^{\delta_{\alpha\beta}+\delta_{\beta\gamma}+\delta_{\gamma\alpha}}\langle
\epsilon_{\alpha}\cdot\Psi,\, \epsilon_{\beta}\cdot
\epsilon_{\gamma}\cdot\Psi\rangle=-\langle \epsilon_{\alpha}\cdot\Psi,\,
\epsilon_{\beta}\cdot \epsilon_{\gamma}\cdot\Psi\rangle
\end{array}
$$
where we have used (10) and (7). It follows that
$$
Re\langle \epsilon_{\alpha}\cdot\Psi,\, \epsilon_{\beta}\cdot
\epsilon_{\gamma}\cdot\Psi\rangle=0.
$$ \hfill $\Box$

\

{\bf Proposition 3.2}\, {\em For a map $\phi:(M,\,g)\to (N,\,h)$ and a spinor
$\Psi\in \Gamma(\Sigma M)$, we define a spinor field
$\psi_{\phi,\,\Psi}$ along the map by {\rm (4)}. Then

{\rm (i)}\,$\mathcal{R}(\phi,\,\psi_{\phi,\,\Psi})\equiv 0$};

{\rm (ii)}\,$\slashed{D}\psi_{\phi,\,\Psi}=-\Psi\otimes
\tau(\phi)-2(\nabla_{\epsilon_{\alpha}}\Psi+\frac
12\epsilon_{\alpha}\cdot\slashed{\partial}\Psi)\otimes\phi_*(\epsilon_{\alpha})$
{\em where $\epsilon_{\alpha}$ ($\alpha=1,\,2$), as always, is a local orthonormal basis
of $M$}.

\

{\bf Remark}\, (a) The Dirac-harmonicity of
$(\phi,\,\psi_{\phi,\,\Psi})$ implies the harmonicity of $\phi$ by
(i) and (8).

(b) $(\nabla_{\epsilon_{\alpha}}\Psi+\frac
12\epsilon_{\alpha}\cdot\slashed{\partial}\Psi)\otimes\phi_*(\epsilon_{\alpha})$
is globally defined.

\

{\bf Proof of Proposition 3.2:}\,\,\,(i) Define local vector fields
$\nabla \phi^i$ on $M$ by
$$
\nabla \phi^i:=(d\phi)^{\sharp}(dy^i)
$$
where $\{dy^i\}$ is the natural local dual basis on $N$. By using
(4), we have
$$
\psi^i:=\psi_{\phi,\,\Psi}(dy^i)=\nabla \phi^i\cdot\Psi
$$
Set
$d\phi=\phi^i_{\alpha}\theta^{\alpha}\otimes\frac{\partial}{\partial
y^i}$ where $\theta^{\alpha}$ is the dual basis for
$\epsilon_{\alpha}$. Then $\nabla \phi^i=\sum
\phi^i_{\alpha}\epsilon_{\alpha}$ and
$$
\langle\psi^k,\,\nabla\phi^j\cdot\psi^l\rangle=\phi^k_{\alpha}\phi^j_{\beta}\phi^l_{\gamma}\langle
\epsilon_{\alpha}\cdot\Psi,\, \epsilon_{\beta}\cdot \epsilon_{\gamma}\cdot\Psi\rangle.
$$
Together with Lemma 3.1, we conclude that
$R^i{}_{jkl}\langle\psi^k,\,\nabla\phi^j\cdot\psi^l\rangle$ is
purely imaginary. On the other hand, from
the proof of Lemma 2.1,
$R^i{}_{jkl}\langle\psi^k,\,\nabla\phi^j\cdot\psi^l\rangle$ must be
real, and hence
$$
\mathcal{R}(\phi,\,\psi_{\phi,\,\Psi})\equiv \frac
12R^i{}_{jkl}\langle\psi^k,\,\nabla\phi^j\cdot\psi^l\rangle\frac{\partial}{\partial
y^i}\equiv 0.
$$

(ii) By using (10) we have
$$
\nabla_{\epsilon_{\alpha}}\Psi+\frac
12\epsilon_{\alpha}\cdot\slashed{\partial}\Psi=\nabla_{\epsilon_{\alpha}}\Psi+\frac
12\epsilon_{\alpha}\cdot\left[\Sigma
\epsilon_{\beta}\cdot\nabla_{\epsilon_{\beta}}\Psi\right]=\left\{\begin{array}{ll}
\frac 12(\nabla_{\epsilon_1}\Psi+\epsilon_{1}\cdot
\epsilon_{2}\cdot\nabla_{\epsilon_2}\Psi),& \alpha=1\\
\frac 12(\nabla_{\epsilon_2}\Psi-\epsilon_{1}\cdot
\epsilon_{2}\cdot\nabla_{\epsilon_1}\Psi),&\alpha=2
\end{array}
\right..\eqno(14)
$$
We choose a local orthonormal frame field $\epsilon_{\alpha}$ such
that $\nabla_{\epsilon_{\alpha}}\epsilon_{\beta}=0$ at $x\in M$.
Then
$$
\begin{array}{ccl}
\slashed{D}\psi_{\phi,\,\Psi}
&=&
\epsilon_{\beta}\cdot\tilde{\nabla}_{\epsilon_{\beta}}\psi_{\phi,\,\Psi}\\
 &=&
\epsilon_{\beta}\cdot\tilde{\nabla}_{\epsilon_{\beta}}
\left(\epsilon_{\alpha}\cdot\Psi\otimes
\phi_*(\epsilon_{\alpha})\right)\\
&=& \epsilon_{\beta}\cdot\left[\nabla_{\epsilon_{\beta}}(\epsilon_{\alpha}\cdot\Psi)\otimes\phi_*(\epsilon_{\alpha})
+\epsilon_{\alpha}\cdot\Psi\otimes\nabla_{\epsilon_{\beta}}(\phi_*(\epsilon_{\alpha}))\right]\\
&=&
\epsilon_{\beta}\cdot\left[((\nabla_{\epsilon_{\beta}}(\epsilon_{\alpha})\cdot\Psi
+\epsilon_{\alpha}\cdot\nabla_{\epsilon_{\beta}}\Psi
 )\otimes\phi_*(\epsilon_{\alpha})
+\epsilon_{\alpha}\cdot\Psi\otimes\nabla_{\epsilon_{\beta}}(\phi_*(\epsilon_{\alpha}))\right]\\
&=&
\epsilon_{\beta}\cdot\epsilon_{\alpha}\cdot\left\{\nabla_{\epsilon_{\beta}}\Psi\otimes\phi_*(\epsilon_{\alpha})
+\Psi\otimes\nabla_{\epsilon_{\beta}}(\phi_*(\epsilon_{\alpha}))\right\}\\
&=&
(\Sigma_{\alpha=\beta}+\Sigma_{\alpha\neq\beta})\epsilon_{\beta}\cdot\epsilon_{\alpha}\cdot\left\{\nabla_{\epsilon_{\beta}}\Psi\otimes\phi_*(\epsilon_{\alpha})
+\Psi\otimes\nabla_{\epsilon_{\beta}}(\phi_*(\epsilon_{\alpha}))\right\}\\
 &=&
 (I)+(II).
\end{array}
\eqno(15)
$$
where
$$
\begin{array}{ccl}
(I) &=&
\epsilon_{\alpha}\cdot\epsilon_{\alpha}\cdot\left\{\nabla_{\epsilon_{\alpha}}\Psi\otimes\phi_*(\epsilon_{\alpha})
+\Psi\otimes\nabla_{\epsilon_{\alpha}}(\phi_*(\epsilon_{\alpha}))\right\}\\
&=&
-\left\{\nabla_{\epsilon_{\alpha}}\Psi\otimes\phi_*(\epsilon_{\alpha})
+\Psi\otimes\left[\nabla_{\epsilon_{\alpha}}(\phi_*(\epsilon_{\alpha}))-\phi_*(\nabla_{\epsilon_{\alpha}}(\phi_*(\epsilon_{\alpha}))\right]\right\}\\
&=&
-\left\{\nabla_{\epsilon_{\alpha}}\Psi\otimes\phi_*(\epsilon_{\alpha})
+\Psi\otimes\tau(\phi)\right\}
\end{array}
\eqno(16)
$$
and
$$
\begin{array}{ccl}
(II) &=&
\epsilon_1\cdot\epsilon_2\cdot\left\{\nabla_{\epsilon_1}\Psi\otimes\phi_*(\epsilon_2)
+\Psi\otimes\nabla_{\epsilon_1}(\phi_*(\epsilon_2))\right\}\\
&&
+\epsilon_2\cdot\epsilon_1\cdot\left\{\nabla_{\epsilon_2}\Psi\otimes\phi_*(\epsilon_1)
+\Psi\otimes\nabla_{\epsilon_2}(\phi_*(\epsilon_1))\right\}\\
&=&
\epsilon_1\cdot\epsilon_2\cdot\left\{\nabla_{\epsilon_1}\Psi\otimes\phi_*(\epsilon_2)-\nabla_{\epsilon_2}\Psi\otimes\phi_*(\epsilon_1)
+\Psi\otimes\nabla_{\epsilon_1}(\phi_*(\epsilon_2))-\Psi\otimes\nabla_{\epsilon_2}(\phi_*(\epsilon_1)\right\}\\
&=&
\epsilon_1\cdot\epsilon_2\cdot\left\{\nabla_{\epsilon_1}\Psi\otimes\phi_*(\epsilon_2)-\nabla_{\epsilon_2}\Psi\otimes\phi_*(\epsilon_1)\right\}
\end{array}
\eqno(17)
$$
here we have used the following
$$
\nabla_{\epsilon_1}(\phi_*(\epsilon_2))=(\nabla_{\epsilon_1}\phi_*)(\epsilon_2)=(\nabla_{\epsilon_2}\phi_*)(\epsilon_1)=\nabla_{\epsilon_2}(\phi_*(\epsilon_1))
$$
Substituting (16) and (17) into (15) yields
$$
\begin{array}{ccl}
\slashed{D}\psi_{\phi,\,\Psi} &=&
-\left\{\nabla_{\epsilon_{\alpha}}\Psi\otimes\phi_*(\epsilon_{\alpha})
+\Psi\otimes\tau(\phi)\right\}+\epsilon_1\cdot\epsilon_2\cdot\left\{\nabla_{\epsilon_1}\Psi\otimes\phi_*(\epsilon_2)-\nabla_{\epsilon_2}\Psi\otimes\phi_*(\epsilon_1)\right\}\\
&=& -\Psi\otimes
\tau(\phi)-(\nabla_{\epsilon_1}\Psi+\epsilon_{1}\cdot
\epsilon_{2}\cdot\nabla_{\epsilon_2}\Psi)\otimes \phi_*(\epsilon_{1})\\
 &&\qquad\qquad\quad+(\epsilon_{1}\cdot
\epsilon_{2}\cdot\nabla_{\epsilon_1}\Psi-\nabla_{\epsilon_2}\Psi)\otimes
\phi_*(\epsilon_{2})
\end{array}
\eqno(18)
$$
Plugging (14) into (18) yields (ii). \hfill $\Box$

\section{Proof of Theorem 2 and Examples}

{\bf Proof of Theorem 2}\,\, By using (i) of Proposition 3.2 and
the harmonicity of $\phi$ we have
$$
\mathcal{R}(\phi,\,\psi_{\phi,\,\Psi})\equiv 0\equiv \tau(\phi).
$$
Thus, $(\phi,\,\psi_{\phi,\,\Psi})$
 satisfies (8). On the other hand, since $\Psi$ is a twistor spinor
 and $n=2$
 we get
 $$
 \nabla_{\epsilon_{\alpha}}\Psi+\frac 12\epsilon_{\alpha}\cdot\slashed{\partial}\Psi=0.
 $$
 Plugging this into the equation in (ii) of Proposition 3.2 yields
 $\slashed{D}\psi_{\phi,\,\Psi}=0$. It follows that $(\phi,\,\psi_{\phi,\,\Psi})$
 satisfies (9), and hence $(\phi,\,\psi_{\phi,\,\Psi})$ is  a
 Dirac-harmonic map. \hfill $\Box$

 \

 {\bf Corollary 4.1}\,\,{\em Let $\psi_{\phi,\,\Psi}$ be defined by
(4) from a branched minimal conformal immersion
$\phi:M\hookrightarrow N$ and a twistor spinor $\Psi\in
\Gamma(\Sigma M)$. Then $(\phi,\,\psi_{\phi,\,\Psi})$ is a
Dirac-harmonic map.}

\

This corollary comes from the fact that a conformal map from a Riemann
surface is harmonic if and only if it is a branched minimal
immersion [6]. Say that an almost Hermitian manifold $(N,\,h,\,J)$ is {\em $(1,\,2)$-symplectic}
if
$$
\nabla^N_{\bar{Z}}W\in \Gamma(T^{1,\,0}N)\quad\mbox{for every}\quad Z,\,W\in \Gamma(T^{1,\,0}N).
$$

 Lichnerowicz proved in [22] that any holomorphic map from a cosymplectic manifold to a
 $(1,\,2)$-symplectic manifold is harmonic. Since a Riemann surface is automatically
 cosymplectic, we have the following:

\

{\bf Corollary 4.2}\,\,{\em Let $\psi_{\phi,\,\Psi}$ be defined by
(4) from a holomorphic map $\phi:M\to N$ and a twistor spinor
$\Psi\in \Gamma(\Sigma M)$ where $N$ is a $(1,\,2)$-symplectic
manifold. Then $(\phi,\,\psi_{\phi,\,\Psi})$ is a Dirac-harmonic
map.}

\

 {\bf Example 1}\,(non-conformal Dirac-harmonic maps)\,\, Suppose that $\mathbb{R}^2$ is given the
 metric $ds^2=2dzd\bar{z}$, where $z=x+iy$ is the standard complex coordinate, and let $\mathbf{e}_0,\,\cdots,\,\mathbf{e}_n$ be a unitary basis of $\mathbb{C}^{n+1}$.
 Define $\phi:\mathbb{R}^2\to \mathbb{C}\mathrm{P}^n$ by
 $$
 \phi(z)=\left[\sum_{j=0}^nr_j\exp(\mu_jz-\overline{\mu_jz})\mathbf{e}_j\right]
 $$
 where $r_0,\,\cdots,\,r_n$ are strictly positive real numbers and $\mu_0,\,\cdots,\,\mu_n$
 are complex numbers of unit modulus satisfying
 $$
 \sum_{j=0}^nr_j^2=1,\quad\qquad \sum_{j=0}^nr_j\mu_j=0.
 $$
 Then $\phi$ is a harmonic map [6,\,19]. In particular, $\phi$ is
 totally real,
and it is conformal if and only if
$$
\sum_{j=0}^nr_j\mu_j^2=0.
$$
Let us consider a twistor spinor $\Psi:\mathbb{R}^2\to \Delta_2=\mathbb{C}^2$  on $\mathbb{R}^2$ (cf [18]).
According to Example 1 of [1] the set of all twistor spinors on $\mathbb{R}^2$ is given by
$$
\Psi(z)=\Psi_0-\frac 12 z\cdot\Psi_1
$$
with $\Psi_0,\,\Psi_1\in \Delta_2$.  From Theorem 2, we obtain that $(\phi,\,\psi_{\phi,\,\Psi})$ is a Dirac-harmonic map
from $\mathbb{R}^2$ into $\mathbb{C}\mathrm{P}^n$
where
$$
\psi_{\phi,\,\Psi}:=\Sigma_{\alpha} \epsilon_{\alpha}\cdot\Psi\otimes
\phi_*(\epsilon_{\alpha})
$$
where $\epsilon_{\alpha}$ ($\alpha=1,\,2$) is a local orthonormal basis
of $M$.
Furthermore, $\phi$ is non-conformal if $\sum_{j=0}^nr_j\mu_j^2\neq 0$.

\

{\bf Example 2}\,\,\,(Dirac-harmonic sequence) For each $p=0,\,\cdots,\, n$,
let $\phi_p:S^2\to \mathbb{C}\mathrm{P}^n$ be given by
$$
\phi_p[z_0,\,z_1]=\left[f_{p,\,0}(z_0/z_1),\,\cdots,f_{p,\,n}(z_0/z_1)\right]
$$
where $[z_0,\,z_1]\in \mathbb{C}\mathrm{P}^1=S^2$, and for $r=0,\, \cdots,\, n$, $f_{p,\,r}(z)$
is given by
$$
f_{p,\,r}(z)=\frac{p\,!}{(1+z\bar{z})^p}\sqrt{C_r^n}z^{r-p}\sum_k(-1)^kC^r_{p-k}C^{n-r}_k(z\bar{z})^k
$$
where
$$
C^n_r=\frac{n(n-1)\cdots(n-r+1)}{r!}.
$$
Then $\phi_p$ is a conformal minimal immersion (therefore it is a harmonic map) with induced metric
$$
ds_p^2=\frac{n+2p(n-p)}{(1+z\bar{z})^2}dzd\bar{z}.
$$
According to Theorem 7 of [1] the twistor spinors on
$(S^2,\,ds_p^2)$ are given by
$$
\Psi(z)=\frac{\Psi_0+z\cdot \Psi_1}{\sqrt{1+z\bar{z}}}
$$
where $\Psi_0,\,\Psi_1\in \Delta_2$ are constants and where we
identify the new and old spin bundles as in [1]. Thus we obtain a
Dirac-harmonic sequence  $(\phi_p,\,\psi_{\phi_p,\,\Psi})$ from
$S^2$ into $\mathbb{C}\mathrm{P}^n$ (cf. [7]) where
$$
\psi_{\phi_p,\,\Psi}:=\Sigma_{\alpha} \epsilon_{\alpha}\cdot\Psi\otimes
\phi_{p*}(\epsilon_{\alpha}).
$$

\section{Dirac-harmonic maps from Riemannian manifolds}

In this section, we are going to construct Dirac-harmonic maps
$(\phi,\,\psi)$ for which $\phi$ is not harmonic.

Let $(N,\,h)$ be a Riemannian manifold of dimension $n'$, $(M,\,g)$ be
an $n$-dimensional Riemannian manifold with fixed spin structure,
$\Sigma M$ its spinor bundle, with induced Hermitian metric $\langle\cdot,\,\cdot\rangle$.
Let $\phi:M\hookrightarrow N$ be an isometric immersion which means that the
natural induced Riemannian metric on $M$ from the ambient space
$N$ coincides with the original one on $M$. We identify $M$ with its immersed image in $N$.
For each $x\in M$ the tangent space $T_xN$ can be decomposed into a direct sum of
$T_xM$ and its orthogonal complement $T_x^{\bot}M$. Such a decomposition is differentiable.
Thus, we have an orthogonal decomposition of the tangent bundle $TN$ along $M$
$$
TN|_M=\phi^{-1}TN=TM\oplus T^{\bot}M.
$$
For a global section $\mathcal{R}(\phi,\,\psi)$ on $\phi^{-1}TN$ (see Section 2), we have
$$
\mathcal{R}(\phi,\,\psi)=\mathcal{R}^T(\phi,\,\psi)+\mathcal{R}^N(\phi,\,\psi)
$$
where
$$
\mathcal{R}^T(\phi,\,\psi)\in \Gamma(TM),\qquad\mathcal{R}^N(\phi,\,\psi)\in \Gamma(T^{\bot}M).
$$
Similarly, for $\slashed{D}\psi\in\Gamma(\Sigma M\otimes\phi^{-1}TN)$, we have
$$
\slashed{D}\psi=\slashed{D}^T\psi+\slashed{D}^N\psi
$$
where
$$
\slashed{D}^T\psi\in\Gamma(\Sigma M\otimes TM),\qquad \slashed{D}^N\psi\in\Gamma(\Sigma M\otimes T^{\bot}M).
$$
The mean curvature vector of $M$ in $N$ is
$$
H=\frac 1n \tau(\phi)\in\Gamma(T^{\bot}M)
$$
where $\tau(\phi)$ is the tension field of the map $\phi$. Hence we have the following:

\

{\bf Lemma 5.1}\,{\em Let $\phi:M\hookrightarrow N$ be an isometric immersion with the mean curvature vector $H$ and $\psi\in\Gamma(\Sigma M\otimes \phi^{-1}TN$). Then $(\phi,\,\psi)$ is a Dirac-harmonic
map from $M$ into $N$ if and only if

(i) $\mathcal{R}^T(\phi,\,\psi)=0$;

(ii) $\mathcal{R}^N(\phi,\,\psi)=nH$ {\em where} $n={\rm dim}\,M$;

(iii) $\slashed{D}^T\psi=0$;

(iv) $\slashed{D}^N\psi=0$.}

\

In this section we shall be using the following ranges of indices:
$$
1\leq\alpha,\,\beta,\,\cdots\leq n,\qquad n+1\leq s,\,t,\,\cdots\leq n',\qquad 1\leq i,\,j,\,\cdots\leq n'.
$$

Choose a local frame field $\{\epsilon_i\}$ of $\phi^{-1}TN$ such
that $\{\epsilon_{\alpha}\}$ lies in the tangent bundle $TM$ and
$\{\epsilon_s\}$ in the normal bundle $T^{\bot}M$ of $M$. By using
(12) we have
$$
\nabla \phi^j=\sum_{\alpha=1}^n
\delta_{\alpha}^j\epsilon_{\alpha}. \eqno(19)
$$
Plugging (19) into (13) yields
$$
\mathcal{R}(\phi,\,\psi)=\frac 12R^i{}_{\alpha
kl}\left(x\right)\langle\psi^k,\,\epsilon_{\alpha}\cdot\psi^l\rangle\epsilon_i\left(x\right).
\eqno(20)
$$
Choose a local orthonormal frame field $\{\epsilon_{\alpha}\}$
near $x\in M$ with
$\nabla_{\epsilon_{\alpha}}\epsilon_{\beta}|_x=0$. By (11) we have
$$
\begin{array}{ccl}
\slashed{D}\psi
&=&
\slashed{D}(\psi^i\otimes\epsilon_i)\\
&=&
\epsilon_{\alpha}\cdot\tilde{\nabla}_{\epsilon_{\alpha}}(\psi^i\otimes\epsilon_i)\\
&=&
\epsilon_{\alpha}\cdot\left[(\nabla_{\epsilon_{\alpha}}\psi^i)\otimes\epsilon_i
+\psi^i\otimes\nabla_{\epsilon_{\alpha}}\epsilon_i\right]\\
&=&
(\epsilon_{\alpha}\cdot\nabla_{\epsilon_{\alpha}}\psi^i)\otimes\epsilon_i
+
\epsilon_{\alpha}\cdot\left[\psi^{\beta}\otimes\nabla_{\epsilon_{\alpha}}\epsilon_{\beta}
+\psi^s\otimes\nabla_{\epsilon_{\alpha}}\epsilon_s\right]\\
&=&
\slashed{\partial}\psi^i\otimes\epsilon_i+\epsilon_{\alpha}\cdot\psi^s\otimes\nabla_{\epsilon_{\alpha}}\epsilon_s
\end{array}
\eqno(21)
$$
at $x$.

Let $A_{\nu}$ be the shape operator and $\nabla^{\bot}_X$ the normal
connection of $M$ in $N$ where $X$ denotes a tangent vector of $M$ and
$\nu$ a normal vector to $M$. Then
$$
\nabla_{\epsilon_{\alpha}}\epsilon_s=-A_{\epsilon_s}\epsilon_{\alpha}+\nabla^{\bot}_{\epsilon_{\alpha}}\epsilon_s.\eqno(22)
$$
Let $B$ be the second fundamental form of $M$ in $N$. Then $B$ satisfies the Weingarten equation
$$
\langle B(X,\,Y),\,\nu\rangle=\langle A_{\nu}(X),\,Y\rangle
\eqno(23)
$$
where $X,\,Y\in \Gamma(TM)$. By using (22) and (23) we have
$$
\nabla_{\epsilon_{\alpha}}\epsilon_s=-\langle
B(\epsilon_{\alpha},\,\epsilon_{\beta}),\,\epsilon_s\rangle\epsilon_{\beta}
+\nabla^{\bot}_{\epsilon_{\alpha}}\epsilon_s.\eqno(24)
$$
By plugging (24) into (21) we obtain
$$
\slashed{D}\psi=\slashed{\partial}\psi^i\otimes\epsilon_i-\langle
B(\epsilon_{\alpha},\,\epsilon_{\beta}),\,\epsilon_s\rangle\epsilon_{\alpha}\cdot\psi^s\otimes\epsilon_{\beta}
+\epsilon_{\alpha}\cdot\psi^s\otimes\nabla^{\bot}_{\epsilon_{\alpha}}\epsilon_s.
\eqno(25)
$$

\

Let $(\cdots)^{T}$ and $(\cdots)^{N}$ denote the orthogonal projection into the tangent bundle
$\Sigma M\otimes TM$ and the normal bundle $\Sigma M\otimes T^{\bot}M$ respectively.

\

{\bf Lemma 5.2}\,\,{\em Let $\psi^{T}$ be defined by
$$
\psi^{T}=\Sigma_{\alpha} \epsilon_{\alpha}\cdot\Psi\otimes
\phi_*(\epsilon_{\alpha})
$$
from an isometric immersion $\phi:M\hookrightarrow N$ and a spinor $\Psi\in
\Gamma(\Sigma M)$ where $\epsilon_{\alpha}$  is a local
orthonormal basis on $M$. Then
$$
\slashed{D}^T\psi=-\left[2\nabla_{\epsilon_{\beta}}\Psi+\epsilon_{\beta}\cdot\slashed{\partial}\Psi
+\langle
B(\epsilon_{\alpha},\,\epsilon_{\beta}),\,\epsilon_s\rangle\epsilon_{\alpha}\cdot\psi^s\right]\otimes\epsilon_{\beta}
\eqno(26)
$$
where $\psi^N=\Sigma_s\psi^s\otimes\epsilon_s$. In particular, if $N=N(c)$ is a Riemannian manifold of constant curvature $c$, then
$$
\mathcal{R}^T(\phi,\,\psi)=0,
$$
$$
\mathcal{R}^N(\phi,\,\psi)=-2nc Re \langle\psi^s,\,\Psi\rangle\epsilon_s
$$
where $n={\rm dim}M$.
}

\

{\bf Proof}\,\, Choose a local orthonormal frame field $\{\epsilon_{\alpha}\}$ near $x\in M$
with $\nabla_{\epsilon_{\alpha}}\epsilon_{\beta}|_x=0$.
$$
\begin{array}{ccl}
\slashed{\partial}\psi^{\alpha}
&=&
\slashed{\partial}(\epsilon_{\alpha}\cdot\Psi)\\
&=&
\epsilon_{\beta}\cdot\nabla_{\epsilon_{\beta}}(\epsilon_{\alpha}\cdot\Psi)\\
&=&
\epsilon_{\beta}\left[(\nabla_{\epsilon_{\beta}}{\epsilon_{\alpha}})\cdot\Psi+
\epsilon_{\alpha}\cdot\nabla_{\epsilon_{\beta}}\Psi\right]\\
&=&
\epsilon_{\beta}\cdot\epsilon_{\alpha}\cdot\nabla_{\epsilon_{\beta}}\Psi\\
&=&
-\nabla_{\epsilon_{\alpha}}\Psi-
\sum_{\beta\neq\alpha}\epsilon_{\alpha}\cdot\epsilon_{\beta}\cdot\nabla_{\epsilon_{\beta}}\Psi\\
&=&
-2\nabla_{\epsilon_{\alpha}}\Psi-\epsilon_{\alpha}\cdot\slashed{\partial}\Psi.
\end{array}
\eqno(27)
$$
Substituting (27) into (25) and taking the tangent projection
yield (26). Now we assume that $N:=N(c)$ is of constant curvature
$c$. Then the components of the Riemannian curvature tensor of $N$
satisfy
$$
R^i{}_{jkl}=c(\delta^i{}_k\delta_{jl}-\delta^i{}_l\delta_{jk}).
$$
From which together with (20) we obtain
$$
\begin{array}{ccl}
\mathcal{R}(\phi,\,\psi)
&=&
c(\delta^i{}_k\delta_{\alpha l}-\delta^i{}_l\delta_{\alpha k})
Re \langle\psi^k,\,\epsilon_{\alpha}\cdot\psi^l\rangle\epsilon_i\\
&=&
c\left[Re \langle\psi^i,\,\epsilon_{\alpha}\cdot\psi^{\alpha}\rangle
-Re \langle\psi^{\alpha},\,\epsilon_{\alpha}\cdot\psi^i\rangle\right]\epsilon_i\\
&=&
2c Re \langle\psi^i,\,\epsilon_{\alpha}\cdot\psi^{\alpha}\rangle \epsilon_i.
\end{array}
$$
It follows that
$$
\begin{array}{ccl}
\mathcal{R}^T(\phi,\,\psi)
&=&
2c Re \langle\psi^{\beta},\,\epsilon_{\alpha}\cdot\psi^{\alpha}\rangle \epsilon_{\beta}\\
&=&
2c Re \langle \epsilon_{\beta}\cdot\Psi,\,\epsilon_{\alpha}\cdot \epsilon_{\alpha}\cdot\Psi\rangle \epsilon_{\beta}\\
&=&
-2c Re \langle \epsilon_{\beta}\cdot\Psi,\,\Psi\rangle \epsilon_{\beta}=0
\end{array}
\eqno(28)
$$
and
$$
\begin{array}{ccl}
\mathcal{R}^N(\phi,\,\psi)
&=&
2c Re \langle\psi^{s},\,\epsilon_{\alpha}\cdot\psi^{\alpha}\rangle \epsilon_{s}\\
&=&
2c Re \langle \psi^s,\,\epsilon_{\alpha}\cdot \epsilon_{\alpha}\cdot\Psi\rangle \epsilon_{s}\\
&=&
-2c Re \langle \psi^s,\,\Psi\rangle \epsilon_{s}.
\end{array}
$$
Here we have used
$$
\overline{\langle \epsilon_{\beta}\cdot\Psi,\, \Psi\rangle}=-\langle \epsilon_{\beta}\cdot\Psi,\, \Psi\rangle.
$$ \hfill $\Box$

\

We call a spinor $\Phi$ {\em harmonic} if it satisfies the Dirac equation without potential [3],
$$
\slashed{\partial}\Phi=0
$$
where $\slashed{\partial}$
is the usual Dirac operator [14].\\

In the rest of this section, we discuss hypersurfaces in a
Riemannian manifold.

\

{\bf Lemma 5.3}\,\,{\em Let $\phi:M\hookrightarrow N$ be an isometric immersion with codimension $1$ and
$\psi\in\Gamma(\Sigma M\otimes \phi^{-1}TN)$ defined by
$$
\psi=\Sigma_{\alpha} \epsilon_{\alpha}\cdot\Psi\otimes
\phi_*(\epsilon_{\alpha})+\Phi\otimes\nu
$$
where $\nu$ is unit normal vector of $M$,  $\Psi,\,\Phi\in \Gamma(\Sigma M)$ and $\epsilon_{\alpha}$  is a local
orthonormal basis of $M$. Then

(i) $$\slashed{D}^T\psi=0$$
if and only if for each $\beta$
$$
2\epsilon_{\beta}\cdot\nabla_{\epsilon_{\beta}}\Psi-\slashed{\partial}\Psi=\lambda_{\beta}\Phi
\eqno(29)
$$
where $\lambda_{\beta}$ is the principal curvature of $M$ in the direction $\epsilon_{\beta}$;

(ii)
$$\slashed{D}^N\psi=0$$
if and only if $\Phi$ is a harmonic
spinor.
}

\

{\bf Proof}\,\, It is easy to see that
$$
\langle B(\epsilon_{\alpha},\,\epsilon_{\beta}),\,\nu\rangle\epsilon_{\alpha}\cdot\Phi\otimes\epsilon_{\beta}
$$
is globally defined. Choose an adapted orthonormal frame of $M$ such that
$$
\langle B(\epsilon_{\alpha},\,\epsilon_{\beta}),\,\nu\rangle=\lambda_{\alpha}\delta_{\alpha\beta}
$$
where $\lambda_{\alpha}$ is the principal curvature of $\phi$.
Plugging this into (26) yields
$$
\slashed{D}^T\psi=-(2\nabla_{\epsilon_{\beta}}\Psi+\epsilon_{\beta}\cdot\slashed{\partial}\Psi+
\lambda_{\beta}\epsilon_{\beta}\cdot \Phi)\otimes\epsilon_{\beta}.
$$
It follow that $\slashed{D}^T\psi=0$ if and only if
$$
2\nabla_{\epsilon_{\beta}}\Psi+\epsilon_{\beta}\cdot\slashed{\partial}\Psi=-
\lambda_{\beta}\epsilon_{\beta}\cdot \Phi \eqno(30)
$$
for each $\beta$. From (10), we see that (30) holds if and only if
(29) holds for each $\beta$.

\

(ii) Note that $M$ is a hypersurface. It follows that
$\nabla^{\bot}\nu=0$. Plugging this into (25) yields
$$
\slashed{D}^N\psi=\slashed{\partial}\Phi\otimes\nu
+\epsilon_{\alpha}\cdot\Phi\otimes\nabla^{\bot}_{\epsilon_{\alpha}}\nu=\slashed{\partial}\Phi\otimes\nu
$$
which immediately implies (ii). \hfill $\Box$

\

{\bf Corollary 5.4}\,\,{\em Let $\phi:M\hookrightarrow N$ be an
isometric immersion with codimension $1$. If $(\phi,\,\psi)$ is a
Dirac-harmonic map then $\Phi$ is a harmonic spinor where
$$
\psi=\Sigma_{\alpha} \epsilon_{\alpha}\cdot\Psi\otimes
\phi_*(\epsilon_{\alpha})+\Phi\otimes\nu
$$
where $\nu$ is unit normal vector of $M$,  $\Psi,\,\Phi\in
\Gamma(\Sigma M)$ and $\epsilon_{\alpha}$  is a local orthonormal
basis of $M$}.

\

{\bf Proof of Theorem 1}\,\,(ii) For a totally umbilical hypersurface $M$, we can assume that
$$
\lambda_1=\lambda_2=\cdots=\lambda_n=\langle H,\,\nu\rangle
\eqno(31)
$$
where $\lambda_{\alpha}$ is the principal curvature of $M$. Note
that $\Psi$ is a twistor spinor. Hence from [1, page 23, Theorem
2] the spinor field $X\cdot\nabla_X\psi$ does not depend on the
unit vector field $X$. Together with (3), we obtain
$$
\epsilon_{1}\cdot\nabla_{\epsilon_{1}}\Psi=\cdots=\epsilon_{n}\cdot\nabla_{\epsilon_{n}}\Psi=\frac
1n\slashed{\partial}\Psi=-\frac{\langle H,\,\nu\rangle}{n-2}\Phi
$$
where $n={\rm dim} M$. It follows that
$$
2\epsilon_{\beta}\cdot\nabla_{\epsilon_{\beta}}\Psi-\slashed{\partial}\Psi
=-\frac{2\langle H,\,\nu\rangle}{n-2}\Phi+\frac{n\langle
H,\,\nu\rangle}{n-2}\Phi=\langle H,\,\nu\rangle\Phi.
$$
Now (ii) can be obtained from (31), Lemma 5.1, Lemma 5.2 and Lemma
5.3 immediately.

(i) For a minimal immersion $\phi$,  we can assume that
$$
\lambda_1=-\lambda_2. \eqno(32)
$$
On the other hand,
$$
2\epsilon_{1}\cdot\nabla_{\epsilon_{1}}\Psi-\slashed{\partial}\Psi=
-[2\epsilon_{2}\cdot\nabla_{\epsilon_{2}}\Psi-\slashed{\partial}\Psi].
$$
Together with (2) and (32) we get (29) for $\beta=1,\,2$. Now (i)
can be obtained from Lemma 5.1, Lemma 5.2 and Lemma 5.3 immediately.
\hfill $\Box$

\

{\bf Example 3}\,\,We consider a totally umbilical hypersurface
$\mathbb{R}^n$ in a hyperbolic space form $\mathbb{H}^{n+1}(-1)$ where
$n\geq 3$. We recall the corresponding construction:
For any two vectors $X$ and $Y$ in $\mathbb{R}^{n+2}$, we set
$$
g(X,\,Y)=\sum_{i=1}^{n+1}X^iY^i-X^{n+2}Y^{n+2}.
$$
We define
$$
\mathbb{H}^{n+1}(-1)=\{x\in \mathbb{R}^{n+2}\,\,|\,\,x_{n+2}>0, g(x,\,x)=-1\}.
$$
Then $\mathbb{H}^{n+1}(-1)$ is a connected simply-connected
hypersurface of $\mathbb{R}^{n+2}$ and the restriction of $g$ to the tangent space of $\mathbb{H}^{n+1}(-1)$
yields a complete Riemannian metric of constant curvature $-1$.

Consider the following small spheres [26]
$$
\mathbb{R}^n:=\{x\in\mathbb{H}^{n+1}(-1)\,|\,\,x_{n+2}=x_{n+1}+1\}.
$$
Then the inclusion map $\phi: \mathbb{R}^n\hookrightarrow
\mathbb{H}^{n+1}(-1)$ is a totally umbilical isometric immersion
with respect to the induced metric. Furthermore its sharp operator
is $A=Id$ [16], that is, its principal curvatures satisfy that
$$
\lambda_1=\cdots=\lambda_n=1.
$$
It follows that $H=\nu$. We take a constant  $\Phi\in \Delta_n$  where
$$
\Delta_n=\mathbb{C}^{2^k}\qquad\mbox{for}\qquad n=2k,\,2k+1
$$
is the vector space of complex $n$ spinors (cf. [14] ). Then $\Phi$
is a harmonic spinor on $\mathbb{R}^n$. Let us consider a twistor
spinor $\Psi:\mathbb{R}^n\to \Delta_n$ on $\mathbb{R}^n$ satisfying
$$
\slashed{\partial}\Psi=-\frac{n}{n-2}\Phi
$$
where $n\geq 3$.
Now we integrate the twistor equation
$$
\begin{array}{ccl}
0&=&\nabla_X\Psi+\frac 1n X\cdot\slashed{\partial}\Psi\\
&=&
\nabla_X\Psi-\frac 1n X\cdot\left(\frac{n}{n-2}\Phi\right)\\
&=&
\nabla_X\Psi+\frac 1{2-n} X\cdot\Phi
\end{array}
$$
along the line $\{sX\,|\, 0\leq s\leq 1\}$, i.e.
$$
\begin{array}{ccl}
\Psi(X)-\Psi(0)
&=&
(\Psi\circ\sigma)(1)-(\Psi\circ\sigma)(0)\\
&=&
\int_0^1\frac{d(\Psi\circ\sigma)}{ds}ds\\
&=&
\int_0^1(\nabla_X\Psi)ds\\
&=&
\int_0^1\frac 1{n-2}X\cdot\Phi ds=\frac 1{n-2}X\cdot\Phi
\end{array}
$$
where $\sigma(s):=sX$ and $\Psi(0)\in \Delta_n$ is constant (cf.
[1]). It is easy to see that the solutions of the equation
$\slashed{\partial}\Psi = - \frac{n}{n-2} \Phi$ are given by
$\Psi(X) = \Psi(0) + \frac{1}{n-2} X\cdot \Phi$ (cf.[1,P29,
Example 1]). Now we will find $\Psi_0:=\Psi(0)$ such that (1)
holds.  Note that $\langle \Phi,\,X\cdot\Phi\rangle$ is purely
imaginary. Hence
$$
\begin{array}{ccl}
\langle \Phi,\,\Psi\rangle
&=&
\langle \Phi,\,\Psi_0+\frac 1{n-2}X\cdot\Phi\rangle\\
&=&
\langle \Phi,\,\Psi_0\rangle+\frac 1{n-2}\langle \Phi,\,X\cdot\Phi\rangle=\langle \Phi,\,\Psi_0\rangle+\frac 1{n-2}Im\langle \Phi,\,X\cdot\Phi\rangle.
\end{array}
$$
It is easy to see that (1) holds when $\Phi,\,\Psi_0\in \Delta_n$ satisfy
$$
Re\langle \Phi,\,\Psi_0\rangle=\frac 12. \eqno(33)
$$
 Thus we obtain that $(\phi,\,\psi)$ is a Dirac-harmonic map from
$\mathbb{R}^n$ into $\mathbb{H}^{n+1}(-1)$ where
$$
\psi(X)=\epsilon_{\alpha}\cdot\left(\Psi_0+\frac 1{n-2}X\cdot\Phi\right)\otimes\phi_*\epsilon_{\alpha}+\Phi\otimes\nu
$$
and $\Phi,\,\Psi_0$ satisfy (33).

\

{\bf Remark}\,\, It is easy to prove that if $\psi^T = \sum
\epsilon_{\alpha}\cdot \Psi \otimes \phi_{*}   ( \epsilon_{\alpha}
)$ and $(\phi, \psi)$ is Dirac-harmonic then $n=2$ implies that
$H=0$. Hence when dim$M=2$, $\Phi=0$, (1) automatically holds, and
(2) holds if and only if $\Psi$ is a twistor spinor.

\

\section{Hypersurfaces with constant principal curvatures in a Riemannian manifold of constant curvature}

In this section, we consider first the following example. Equipped
with the pseudo-Riemannian metric
$$
ds^2=dx^2_1+\cdots+dx^2_{n+1}-dx^2_{n+2},
$$
$\mathbb{R}^{n+2}$ becomes {\em Minkowski space}
$\mathbb{R}^{n+2}_1$. We define {\em (real) hyperbolic space}
$$
H^{n+1}(R):=\left\{x\in \mathbb{R}^{n+2}\,|\,q(x)=-R^2, \,
x_{n+2}>0\right\}
$$
where $q(x):=x_1^2+\cdots +x_{n+1}^2-x_{n+2}^2$. $H^{n+1}(R)$ is a
connected, simply-connected hypersurface of $\mathbb{R}^{n+2}_1$
and  the restriction of $ds^2$ to tangent vectors yields a
(positive-definite) complete Riemannian metric in $H^{n+1}(R)$ of
constant sectional curvature $c=-\frac 1{R^2}$.
 We now define a family of product hypersurfaces
$$
M:=\left\{x\in
H^{n+1}(R)\,|\,x^2_1+\cdots+x^2_{k+1}=r^2\right\}=S^k(r)\times
H^{n-k}(\sqrt{R^2+r^2})\eqno(34)
$$
for $r>0$ and $k=1,\cdots,n-1$.  The induced metric on $M$ is
$$
ds^2_{S^k(r)}+ds^2_{H^{n-k}(\sqrt{R^2+r^2})} = r^2
ds^2_{S^k(1)}+(R^2+r^2)ds^2_{H^{n-k}(1)}. \eqno(35)
$$
$M$ has principal curvatures $\frac{\sqrt{R^2+r^2}}{rR}$ with
multiplicity $k$ and $\frac{r}{R\sqrt{R^2+r^2}}$ with multiplicity
$n-k$ [25]. Therefore, the trace of the shape operator of $M$ in
$H^{n+1}(R)$ is $\frac{kR^2+nr^2}{Rr\sqrt{R^2+r^2}}$. We have the
following:

\

{\bf Lemma 6.1}\, {\em Let $M:=S^k(r)\times
H^{n-k}(\sqrt{R^2+r^2})$ be a hypersurface in $H^{n+1}(R)\subset
\mathbb{R}^{n+2}_1$. Then

(i) $M$ is non-minimal, therefore, $\phi:M\looparrowright
H^{n+1}(R)$ is not harmonic;

(ii) $M$ has two constant principal curvatures, with constant
multiplicities.}

\

In order to getting new non-trivial Dirac-harmonic maps, we
construct a spinor field along a hypersurface with two constant
principal curvatures in a Riemannian manifold of constant
curvature. We shall be using the following ranges of indices:
$$
1\leq i,\,j,\,\cdots\leq k,\qquad k+1\leq r,\,s,\,\cdots\leq
n,\qquad 1\leq \alpha,\,\beta,\,\cdots\leq n.
$$

\

{\bf Lemma 6.2}\, {\em Let $\phi:M^n\looparrowright N^{n+1}(c)$ be
a hypersurface with two principal curvatures $\lambda$ and $\mu$
in a Riemannian manifold of constant curvature $c$, where
$\lambda$ has the multiplicity $k$ and $\mu$ has the multiplicity
$n-k$, and $\psi\in\Gamma(\Sigma M\otimes \phi^{-1}TN)$ is defined
by
$$
\psi=\Sigma_{i} \epsilon_{i}\cdot\Psi\otimes
\phi_*(\epsilon_{i})+\Sigma_{r} \epsilon_{r}\cdot\Phi\otimes
\phi_*(\epsilon_{r})+\chi\otimes\nu
$$
where $\nu$ is the unit normal vector of $M$,
$\Psi,\,\Phi,\,\chi\in \Gamma(\Sigma M)$ and $\epsilon_{\alpha}$
is a local orthonormal basis of $M$ such that $\epsilon_i$ is the
eigenvector of $\lambda$ and $\epsilon_r$ is the eigenvector of
$\mu$. Then
$$
\mathcal{R}^T(\phi,\,\psi)
=2c\left[Re\langle\epsilon_i\cdot\Phi,\,\Psi\rangle\epsilon_i-Re\langle\epsilon_r\cdot\Phi,\,\Psi\rangle\epsilon_r\right];
\eqno(36)
$$
$$
\mathcal{R}^N(\phi,\,\psi)=-2c
Re\langle\chi,\,k\Psi+(n-k)\Phi\rangle\nu; \eqno(37)
$$
$$
\slashed{D}^T\psi=-(2\nabla_{\epsilon_{i}}\Psi+\epsilon_{i}\cdot\slashed{\partial}\Psi+
\lambda_{\beta}\epsilon_{i}\cdot \chi)\otimes\epsilon_{i}
-(2\nabla_{\epsilon_{r}}\Phi+\epsilon_{r}\cdot\slashed{\partial}\Phi+
\mu_{\beta}\epsilon_{r}\cdot \chi)\otimes\epsilon_{r}; \eqno(38)
$$
$$
\slashed{D}^N\psi=(\slashed{\partial}\chi)\otimes\nu. \eqno(39)
$$
}

{\bf Proof:}\,\,\, Denote the distributions of the spaces of
principal vectors corresponding to $\lambda$ and $\mu$ by
$D_{\lambda}$ and $D_{\mu}$, i.e.
$$
D_{\lambda}:=\left\{v\in TM\,|\, Av=\lambda v\right\},\qquad
D_{\mu}:=\left\{v\in TM\,|\, Av=\mu v\right\}
$$
where $A$ is the shape operator of $\phi$. Then
$$
\epsilon_i\in D_{\lambda},\qquad \epsilon_r\in D_{\mu} \eqno(40)
$$
and $\psi$ is well-defined. Note that the multiplicities of the
two principal curvatures are constant. Thus $D_{\lambda}$ and
$D_{\mu}$ are completely integrable [24]. In particular, we may
choose a local orthonormal frame field
$\{\epsilon_{i},\,\epsilon_{r}\}$ near $x$ with
$\nabla_{\epsilon_{\alpha}}\epsilon_{\beta}|_x=0$ and satisfying
(40).

Denote $\psi^{T}$ by
$$
\psi^{T}=\psi^{\alpha}\otimes \phi_*(\epsilon_{\alpha}).
$$
Then
$$
\begin{array}{ccl}
\slashed{\partial}\psi^{i} &=&
\slashed{\partial}(\epsilon_{i}\cdot\Psi)\\
&=&
\epsilon_{\beta}\cdot\nabla_{\epsilon_{\beta}}(\epsilon_{i}\cdot\Psi)\\
&=&
\epsilon_{\beta}\left[(\nabla_{\epsilon_{\beta}}{\epsilon_{i}})\cdot\Psi+
\epsilon_{i}\cdot\nabla_{\epsilon_{\beta}}\Psi\right]\\
&=&
\epsilon_{\beta}\cdot\epsilon_{i}\cdot\nabla_{\epsilon_{\beta}}\Psi\\
&=& -\nabla_{\epsilon_{i}}\Psi-
\sum_{\beta\neq i}\epsilon_{i}\cdot\epsilon_{\beta}\cdot\nabla_{\epsilon_{\beta}}\Psi\\
&=&
-2\nabla_{\epsilon_{i}}\Psi-\epsilon_{i}\cdot\slashed{\partial}\Psi.
\end{array}
\eqno(41)
$$
Similarly we have
$$
\slashed{\partial}\psi^{r}=-2\nabla_{\epsilon_{r}}\Phi-\epsilon_{r}\cdot\slashed{\partial}\Phi.
\eqno(42)
$$
By using (25), we have
$$
\slashed{D}\psi=\slashed{\partial}\psi^{\alpha}\otimes\epsilon_{\alpha}+(\slashed{\partial}\chi)\otimes\nu-\langle
B(\epsilon_{\alpha},\,\epsilon_{\beta}),\,\nu\rangle\epsilon_{\alpha}\cdot\chi\otimes\epsilon_{\beta}.
\eqno(43)
$$
Here we have used $\nabla^{\bot}\nu=0$. Plugging (41) and (42)
into (43) and using the Weingarten equation yield
$$
\begin{array}{ccl}
\slashed{D}\psi &=&
-\left[2\nabla_{\epsilon_{i}}\Psi+\epsilon_{i}\cdot\slashed{\partial}\Psi+B(\epsilon_{\alpha},\,\epsilon_{i}),\,\nu\rangle\epsilon_{\alpha}\cdot\chi\right]
\otimes\epsilon_{i}\\
&&-\left[2\nabla_{\epsilon_{r}}\Phi+\epsilon_{r}\cdot\slashed{\partial}\Phi+B(\epsilon_{\alpha},\,\epsilon_{r}),\,\nu\rangle\epsilon_{\alpha}\cdot\chi\right]
\otimes\epsilon_{r}+(\slashed{\partial}\chi)\otimes\nu\\
&=&
-(2\nabla_{\epsilon_{i}}\Psi+\epsilon_{i}\cdot\slashed{\partial}\Psi+
\lambda_{\beta}\epsilon_{i}\cdot \chi)\otimes\epsilon_{i}\\
&&
-(2\nabla_{\epsilon_{r}}\Phi+\epsilon_{r}\cdot\slashed{\partial}\Phi+
\mu_{\beta}\epsilon_{r}\cdot
\chi)\otimes\epsilon_{r}+(\slashed{\partial}\chi)\otimes\nu.
\end{array}
$$
Thus we obtain (38) and (39).

Note that $N^{n+1}(c)$ has constant sectional curvature $c$.
Consider $\epsilon_{\alpha},\,\nu$ as a local orthonormal frame
field of $\phi^{-1}TN$. By simple calculations, we have
$$
\mathcal{R}^T(\phi,\,\psi) =2c
Re\langle\psi^{\beta},\,\epsilon_{\alpha}\cdot\psi^{\alpha}\rangle
\epsilon_{\beta}, \eqno(44)
$$
$$
\mathcal{R}^N(\phi,\,\psi)= 2c Re
\langle\chi,\,\epsilon_{\alpha}\cdot\psi^{\alpha}\rangle \nu.
\eqno(45)
$$
By using the skew-symmetry relation of the Clifford product and
the property of Hermitian metric we have
$$
Re\langle\psi^{i},\,\epsilon_{j}\cdot\psi^{j}\rangle =Re \langle
\epsilon_{i}\cdot\Psi,\,\epsilon_{j}\cdot
\epsilon_{j}\cdot\Psi\rangle =-Re \langle
\epsilon_{i}\cdot\Psi,\,\Psi\rangle=0, \eqno(46)
$$
$$
\begin{array}{ccl}
Re\langle\psi^{i},\,\epsilon_{r}\cdot\psi^{r}\rangle &=&
-Re \langle \epsilon_{i}\cdot\Psi,\,\Phi\rangle\\
&=&
Re \langle \Psi,\,\epsilon_{i}\cdot\Phi\rangle\\
&=& Re\overline{\langle \epsilon_{i}\cdot\Phi,\,
\Psi\rangle}=Re\langle \epsilon_{i}\cdot\Phi,\, \Psi\rangle.
\end{array}
\eqno(47)
$$
Similarly, we have
$$
Re\langle\psi^{r},\,\epsilon_{i}\cdot\psi^{i}\rangle=-Re\langle
\epsilon_{r}\cdot\Phi,\, \Psi\rangle, \eqno(48)
$$
$$
Re\langle\psi^{r},\,\epsilon_{s}\cdot\psi^{s}\rangle=0. \eqno(49)
$$
Substituting (46), (47), (48) and (49) into (44) yields
$$
\begin{array}{ccl}
\mathcal{R}^T(\phi,\,\psi) &=& 2c
Re\langle\psi^{i},\,\epsilon_{j}\cdot\psi^{j}\rangle\epsilon_i
+2c Re\langle\psi^{i},\,\epsilon_{r}\cdot\psi^{r}\rangle\epsilon_i\\
&& +2c
Re\langle\psi^{r},\,\epsilon_{i}\cdot\psi^{i}\rangle\epsilon_r
+2c Re\langle\psi^{r},\,\epsilon_{s}\cdot\psi^{s}\rangle\epsilon_r\\
&=& 2c\left(Re\langle \epsilon_{i}\cdot\Phi,\,
\Psi\rangle\epsilon_i-Re\langle \epsilon_{r}\cdot\Phi,\,
\Psi\rangle\epsilon_r\right).
\end{array}
$$
Finally, using (10) and (45) we obtain (37). \hfill $\Box$

\

\section{Dirac-harmonic maps from surfaces\,\,II}

In this section, we give first a useful criterion for a class of
maps from surfaces into a three-dimensional Riemannian manifold of
constant curvature to be Dirac-harmonic. By using this criterion we
manufacture Dirac-harmonic maps $(\phi,\,\psi)$ from surfaces for
which $\phi$ is not harmonic.

\

{\bf Theorem 7.1}\,{\em Let $\phi:M^2\looparrowright N^{3}(c)$ be
a surface with two principal curvatures $\lambda$ and $\mu$ in a
Riemannian manifold of constant curvature $c$, where $\lambda\neq
\mu$, and let $\psi\in\Gamma(\Sigma M\otimes \phi^{-1}TN)$ be
defined by
$$
\psi= \epsilon_{1}\cdot\Psi\otimes
\phi_*(\epsilon_{1})-\frac{\mu}{\lambda}\epsilon_{2}\cdot\Psi\otimes
\phi_*(\epsilon_{2})+\chi\otimes\nu
$$
where $\nu$ is unit normal vector of $M$,  $\Psi,\,\chi\in
\Gamma(\Sigma M)$ and $\{\epsilon_1,\,\epsilon_2\}$  is a local
orthonormal basis of $M$ such that $\epsilon_1$ is the eigenvector
of $\lambda$ and $\epsilon_2$ is the eigenvector of $\mu$. Then
$(\phi,\,\psi)$ is a Dirac-harmonic map from $M$ into $N$ if and
only if

(i) $\chi$ is a harmonic spinor;

(ii) $c(\frac{\mu}{\lambda}-1) Re\langle\chi,\,\Psi\rangle\nu=H$;

(iii)
$\epsilon_{1}\cdot\nabla_{\epsilon_{1}}\Psi-\epsilon_{2}\cdot\nabla_{\epsilon_{2}}\Psi=\lambda\chi$.}

\

{\bf Proof}\,\, Take $\Phi=-\frac{\mu}{\lambda}\Psi$. Substituting
this into (36) we get
$$
\begin{array}{ccl}
\mathcal{R}^T(\phi,\,\psi) &=&
2c\left[Re\langle\epsilon_1\cdot(-\frac{\mu}{\lambda}\Psi),\,\Psi\rangle\epsilon_1
-Re\langle\epsilon_2\cdot(-\frac{\mu}{\lambda}\Psi),\,\Psi\rangle\epsilon_2\right]\\
&=&
2c\frac{\mu}{\lambda}\left[Re\langle\epsilon_2\cdot\Psi,\,\Psi\rangle\epsilon_2
-Re\langle\epsilon_1\cdot\Psi,\,\Psi\rangle\epsilon_1\right]=0.
\end{array}
$$
Let us assume that (i) (ii) and (iii) hold. From (37) we have
$$
\mathcal{R}^N(\phi,\,\psi)=-2c
Re\langle\chi,\,\Psi-\frac{\mu}{\lambda}\Psi\rangle\nu
=2c(\frac{\mu}{\lambda}-1) Re\langle\chi,\,\Psi\rangle\nu=2H.
$$
By using (39) we obtain
$$
\slashed{D}^N\psi=(\slashed{\partial}\chi)\otimes\nu=0.
$$
From (iii) we have
$$
\epsilon_{2}\cdot\nabla_{\epsilon_{2}}(-\frac{\mu}{\lambda}\Psi)-
\epsilon_{1}\cdot\nabla_{\epsilon_{1}}(-\frac{\mu}{\lambda}\Psi)
=\frac{\mu}{\lambda}\left[\epsilon_{1}\cdot\nabla_{\epsilon_{1}}\Psi-\epsilon_{2}\cdot\nabla_{\epsilon_{2}}\Psi\right]
=\mu\chi.
$$
Together with (38) and (iii) we have
$$
\begin{array}{ccl}
\slashed{D}^T\psi &=&
-(2\nabla_{\epsilon_{1}}\Psi+\epsilon_{1}\cdot\slashed{\partial}\Psi+
\lambda\epsilon_{1}\cdot \chi)\otimes\epsilon_{1}\\
&&
-\left(2\nabla_{\epsilon_{2}}(-\frac{\mu}{\lambda}\Psi)+\epsilon_{2}\cdot\slashed{\partial}(-\frac{\mu}{\lambda}\Psi)+
\mu\epsilon_{2}\cdot \chi\right)\otimes\epsilon_{2}\\
&=&
\epsilon_1\cdot(2\epsilon_1\cdot\nabla_{\epsilon_{1}}\Psi-\slashed{\partial}\Psi-
\lambda\chi)\otimes\epsilon_{1}\\
&&
+\epsilon_2\cdot\left(2\epsilon_2\nabla_{\epsilon_{2}}(-\frac{\mu}{\lambda}\Psi)-\slashed{\partial}
(-\frac{\mu}{\lambda}\Psi)-
\mu\chi\right)\otimes\epsilon_{2}\\
&=&
\epsilon_1\cdot(\epsilon_1\cdot\nabla_{\epsilon_{1}}\Psi-\epsilon_2\cdot\nabla_{\epsilon_{2}}\Psi-
\lambda\chi)\otimes\epsilon_{1}\\
&&
+\epsilon_2\cdot\left(\epsilon_2\cdot\nabla_{\epsilon_{2}}(-\frac{\mu}{\lambda}\Psi)
-\epsilon_1\cdot\nabla_{\epsilon_{1}}(-\frac{\mu}{\lambda}\Psi)-
\mu\chi\right)\otimes\epsilon_{2}\\
&=& (\epsilon_1\cdot
0)\otimes\epsilon_1+\frac{\mu}{\lambda}(\epsilon_2\cdot
0)\otimes\epsilon_2=0.
\end{array}
$$
From Lemma 5.1 we see that $(\phi,\,\psi)$ is a Dirac-harmonic
map.

Conversely, if $(\phi,\,\psi)$ is a Dirac-harmonic map, then it is
easy to see from Lemma 5.1 that (i) (ii) and (iii) hold. \hfill
$\Box$

\

{\bf Remark}\,\, In fact $\phi: S^{1}( r ) \times
H^{1}(\sqrt{R^2+r^2}) \hookrightarrow   H^{3}( R )$ has two
constant principal curvatures $\lambda$ and $\mu$, where $\lambda
\neq \mu$ ( see proof of Theorem 3 below).

\

{\bf Proof of Theorem 3} Let
$$
M:=S^1(r)\times H^1(\sqrt{R^2+r^2})=\{(x_1,\,y_1;\,x_2,\,y_2)\,|\,
x_1^2+y_1^2=r^2,\,x_2^2-y_2^2=-(R^2+r^2),\,y_2>0\}
$$
be parameterized as
$$
\begin{array}{l}
M=(\mathbb{R}/2\pi r\mathbb{Z})\times \mathbb{R}=\\
\left\{(r\cos\frac{\theta}r,\,r\sin\frac{\theta}r,\,\sqrt{R^2+r^2}\sinh\frac
t{\sqrt{R^2+r^2}},\,\sqrt{R^2+r^2}\cosh\frac
t{\sqrt{R^2+r^2}})\,|\,(\theta,\,t)\in [0,\,2\pi r)\times
\mathbb{R}\right\}.
\end{array}
\eqno(50)
$$
The induced metric on $M$ is the  flat metric
$$
d\theta^2+dt^2. \eqno(51)
$$

Since $2^{\left[\frac{dim M}2\right]}=2$, we use "two-component"
spinors. We identify the "untwisted" spinor bundle on $M$ with
$\left[(\mathbb{R}/2\pi r\mathbb{Z})\times
\mathbb{R}\right]\times\mathbb{C}^2$, that is to say, the spinor
on $M$ is a single periodic spinor on $\mathbb{R}^2$ [1, 21]. Let
$\epsilon_1=\frac{\partial}{\partial\theta}$ and
$\epsilon_2=\frac{\partial}{\partial t}$. Then $\epsilon_1$ and
$\epsilon_2$ acting on spinor fields can be identified by
multiplication with matrices [9, 10]
$$
\sigma_1=\left(\begin{array}{ll}  0&1\\
                                      -1&0
                                      \end{array}
                                      \right),\qquad\sigma_2=\left(\begin{array}{ll}  0&i\\
                                      i&0
                                      \end{array}
                                      \right),\qquad i=\sqrt{-1}.
                                      \eqno(52)
$$
Since the metric (51) is flat, $\nabla=d$ is the Levi-Civita
connection on $1$-forms. Hence we have
$$
\nabla_{\epsilon_{\alpha}}=\epsilon_{\alpha} \eqno(53)
$$
on $\Sigma M$ [17, 22].

We take a constant $\chi
=\left(\begin{array}{l}  a\\
b
\end{array}
\right)\in \mathbb{C}^2$. Then $\chi$ is a harmonic spinor on $M$.
Let us consider a spinor
$$
\Psi=\left(\begin{array}{l}  f\\
g
\end{array}
\right):(\mathbb{R}/2\pi r\mathbb{Z})\times \mathbb{R}\to
\mathbb{C}^2 \eqno(54)
$$
on $M$ satisfying
$$
\epsilon_{1}\cdot\nabla_{\epsilon_{1}}\Psi-\epsilon_{2}\cdot\nabla_{\epsilon_{2}}\Psi=\lambda\chi
\eqno(55)
$$
where $\lambda$ is a real constant. By using (52) and (53), (55)
is equivalent to
$$
\left(\frac{\partial}{\partial\theta}-i\frac{\partial}{\partial
t}\right)g=\lambda a, \eqno(56)
$$
$$
-\left(\frac{\partial}{\partial\theta}+i\frac{\partial}{\partial
t}\right)f=\lambda b. \eqno(57)
$$
For arbitrary $m\in \{0,\,1,\,\cdots\}$, we  consider
$g:(\mathbb{R}/2\pi r\mathbb{Z})\times \mathbb{R}\to \mathbb{C}$
defined by
$$
g(\theta,\,t)=\sum_{k=-m}^m e^{i\frac kr\theta}g_k(t). \eqno(58)
$$
From (56) and (58) we have
$$
\lambda a=\frac{\partial g}{\partial\theta}-i\frac{\partial
g}{\partial t} =\frac ir \sum_{k=-m}^m ke^{i\frac
kr\theta}g_k(t)-i\sum_{k=-m}^m e^{i\frac kr\theta}g'_k(t) =\frac
ir \sum_{k=-m}^m e^{i\frac kr\theta}\left[kg_k(t)-rg'_k(t)\right].
$$
Hence we take $g_k$ satisfying
$$
\left\{\begin{array}{ll}
-ig'_k(t)=\lambda a&\quad\mbox{for}\quad k=0\\
kg_k(t)-rg'_k(t)=0&\quad\mbox{for}\quad k\neq 0
\end{array}
\right.. \eqno(59)
$$
One can verify that
$$
g_k(t):=\left\{\begin{array}{ll}
i\lambda at+c_0&\quad\mbox{for}\quad k=0\\
c_ke^{\frac kr t}&\quad\mbox{for}\quad k\neq 0
\end{array}
\right.
$$
satisfies (59). Plugging this into (58) yields
$$
g(\theta,\,t)=\sum_{k=-m}^m c_ke^{\frac kr(t+i\theta)}+i\lambda
at. \eqno(60)
$$
Similarly, the following function
$$
f:(\mathbb{R}/2\pi r\mathbb{Z})\times \mathbb{R}\to \mathbb{C}
$$
defined by
$$
f(\theta,\,t)=\sum_{k=-m}^m d_ke^{\frac kr(-t+i\theta)}+i\lambda
bt \eqno(61)
$$
satisfies (57). Plugging (60) and (61) into (54) yields
$$
\Psi=i\lambda t\left(\begin{array}{l}  b\\
a
\end{array}
\right) +\sum_{k=-m}^m e^{i\frac kr\theta}
\left(\begin{array}{l}  d_ke^{-\frac krt}\\
c_ke^{\frac krt}
\end{array}
\right).
$$
Consider $\phi:M=S^1(r)\times H^{1}(\sqrt{R^2+r^2})\hookrightarrow
H^{3}(R)$. Then $H^{3}(R)$ has constant sectional curvature
$c=-\frac 1{R^2}$. The principal curvatures of $\phi$ are (cf.
Section 6)
$$
\lambda=\frac{\sqrt{R^2+r^2}}{rR},\qquad\mu=\frac{r}{R\sqrt{R^2+r^2}}
$$
and therefore the mean curvature of $\phi$ is
$$
\xi=\frac{R^2+2r^2}{2Rr\sqrt{R^2+r^2}}.
$$
By a straightforward computation one obtains
$$
\left[c\left(\frac{\lambda}{\mu}-1\right)\right]^{-1}\xi=\frac{\sqrt{R^2+r^2}(R^2+2r^2)}{2rR}.
$$

Now we will find $c_k$ and $d_k$ such that (ii) in Theorem 7.1
holds.
$$
\langle\chi,\,\Psi\rangle=a\overline{\left(i\lambda
bt+\sum_{k=-m}^m d_ke^{\frac kr(-t+i\theta)}\right)} +b
\overline{\left(i\lambda at+\sum_{k=-m}^m c_ke^{\frac
kr(t+i\theta)}\right)} =(I) +(II)
$$
where
$$
(I)=a\overline{i\lambda bt}+b\overline{i\lambda
at}=-i(a\bar{b}+b\bar{a})\lambda=-2\lambda iRe(a\bar{b}),
$$
$$
(II)=a\sum_{k=-m}^m\overline{d_ke^{\frac
kr(-t+i\theta)}}+b\sum_{k=-m}^m\overline{c_ke^{\frac
kr(t+i\theta)}}.
$$
Note that (I) is purely imaginary. Hence
$$
\begin{array}{ccl}
Re\langle\chi,\,\Psi\rangle &=&
Re(II)\\
&=&\sum_{k=-m}^m Re\left[a\bar{d}_ke^{-\frac kr(t+i\theta)}\right]
+\sum_{k=-m}^m Re\left[b\bar{c}_ke^{\frac kr(t-i\theta)}\right]\\
&=&\sum_{k=-m}^m Re\left[a\bar{d}_ke^{-\frac kr(t+i\theta)}\right]
+\sum_{k=-m}^m Re\left[b\bar{c}_{-k}e^{-\frac kr(t-i\theta)}\right]\\
&=&
\sum_{k=-m}^m Re\left[a\bar{d}_ke^{-i\frac kr\theta}+b\bar{c}_{-k}e^{i\frac kr\theta}\right]e^{-\frac kr t}\\
&=& \sum_{k=-m}^m Re\left[(a\bar{d}_k+\bar{b}{c}_{-k})e^{-i\frac
kr\theta}\right]e^{-\frac kr t}.
\end{array}
$$
It follows that the sufficient conditions on $c_k$ and $d_k$ for
(ii) in Theorem 7.1 to hold are
$$
Re(a\bar{d_0}+\bar{b}{c_0})=\frac{\sqrt{R^2+r^2}(R^2+2r^2)}{2rR}
$$
and
$$
a\bar{d_k}+\bar{b}c_{-k}=0
$$
for $k=\pm 1,\cdots, \pm m$. \hfill $\Box$

\

{\em Acknowledgements}: The second author would like to thank the Max Planck Institute for Mathematics in the Sciences,
Leipzig, for its hospitality during the preparation of this paper.

\end{document}